\documentclass[]{amsart}

\usepackage{latexsym}
\usepackage{amssymb,amsmath,amsopn}
\usepackage[dvips]{graphicx}   
\usepackage{color,epsfig}      
\usepackage{url}
\usepackage{amscd}
\usepackage{color}

\newtheorem{theorem}{Theorem}[]

\newtheorem{defi}{Definition}[]

\theoremstyle{definition}

\newcommand{\R}{\mathrm{R}}
\newcommand{\Q}{\mathrm{Q}}
\newcommand{\E}{\mathrm{E}}

\newcommand{\Bsp}{\boldsymbol{p}}
\newcommand{\bmf}{\mathbf{f}}

\renewcommand{\Re}{\mathrm{ Re}}

\title{Constructive curves in non-Euclidean planes }
\author[\'A. G.Horv\'ath]{\'Akos G.Horv\'ath}
\address {Department of Geometry \\
Budapest University of Technology and Economics\\
H-1521 Budapest\\
Hungary}
\email{ghorvath@math.bme.hu}

\date{Sept, 2016}

\begin{document}

\begin{abstract}
In this paper we overview the theory of conics and roulettes in four non-Euclidean planes. We collect the literature about these classical concepts, from the eighteenth century to the present, including papers available only on arXiv. The comparison of the four non-Euclidean planes, in terms of the known results on conics and roulettes, reflects only the very subjective view of the author.
\end{abstract}

\maketitle

\emph{Keywords:} Cayley-Klein geometries; conics; Euler-Savary equations; normed plane; roulettes

\emph{2010 MSC:} 46B20, 51M05, 52A21, 53A17

\section{Introduction}

Our purpose is to compare the basic properties of certain non-Euclidean planes through some known results on two types of curves, which is the motivation for the long introduction before the essential part of this paper. These planes are the so-called hyperbolic, spherical, Minkowski and Lorentzian planes. We note that even though there is an analytic approach for the examined curves in most of the planes, which we investigate here, the geometric method of construction leads to the definitions that can be interpreted in every case.

We do not think that everybody regards our approach of these non-Euclidean planes good. We think even less that the collected literature in this paper is complete. However, the large number of rediscovered results in this evergreen topic shows that our paper may not be useless.

In the paper, we first present a unified model of the plane geometries, then compare the theories of the examined constructive curves. Clearly, there is no room here to mention all results discovered in the last three hundred years. Nevertheless, in our approach some important characteristics are shown which behave analogously in all these geometries, and we concentrate on them. In our investigation, it will be very important to know which definitions work simultaneously in all cases, how we can classify the conics, and how basic kinematics works in these planes. We add a large number of references to the introduction of the theory of the curves.

We require an analytic approach to the investigated geometries, and thus, we cannot mention here non-Euclidean planes with only a synthetic model, for example non-Paschian planes based on a synthetic construction as in \cite{gho 8}. The common root of the four geometries used in the paper is a $3$-dimensional vector space endowed with a general scalar multiplication of the vectors. This multiplication is a function with two variables in $V$, which, in the case that $V$ is a real vector space, we call a \emph{product}. First we review a way to axiomatize it. This method gives a natural matching of these geometries, the Lorentzian and Minkowski plane are affine planes, respectively and the hyperbolic and spherical planes are hypersurfaces of certain affine spaces, respectively. The following two subsections contains the comparing of the planes in this direction.

 \subsection{The two Minkowski planes as affine spaces and their common building up}

A generalization of the inner product and the inner product spaces (briefly i.p spaces) was raised by G. Lumer in \cite{lumer}.
The \emph{semi-inner-product (s.i.p)} on a complex vector space $V$ is a complex function $[x,y]:V\times V\longrightarrow \mathbb{C}$ with the following properties:
\begin{description}
\item[s1]: $[x+y,z]=[x,z]+[y,z]$,
\item[s2]: $[\lambda x,y]=\lambda[x,y]$ \mbox{ for every } $\lambda \in \mathbb{C}$,
\item[s3]: $[x,x]>0$ \mbox{ when } $x\not =0$,
\item[s4]: $|[x,y]|^2\leq [x,x][y,y]$,
\end{description}
A vector space $V$ with a s.i.p. is an \emph{s.i.p. space}.

G. Lumer proved that an s.i.p space is a normed vector space with norm $\|x\|=\sqrt{[x,x]}$ and, on the other hand, that every normed vector space can be represented as an s.i.p. space. In \cite{giles} J. R. Giles showed that the following homogeneity property holds:

\begin{description}

\item[s5]: $[x,\lambda y]=\bar{\lambda}[x,y]$ for all complex $\lambda $.

\end{description}
This can be imposed, and all normed vector spaces can be represented as s.i.p. spaces with this property. Giles also introduced the concept of \emph{ continuous s.i.p. space} as an s.i.p. space having the additional property

\begin{description}
\item[s6]: For any unit vectors $x,y \in S$, $\Re\{[y,x+\lambda y]\}\rightarrow\Re\{[y,x]\}$ for all real $\lambda\rightarrow 0$.
\end{description}
The space is uniformly continuous if the above limit is reached uniformly for all points $x,y$ of the unit sphere $S$.

A characterization of the continuous s.i.p. space is based on the differentiability property of the space (see in \cite{giles}).
A normed space is G\^{a}teaux differentiable if for all elements $x,y$ of its unit sphere and real values $\lambda$, the limit
$$
\lim\limits_{\lambda \rightarrow 0}\frac{\|x+\lambda y\|-\|x\|}{\lambda}
$$
exists. A normed vector space is uniformly Fr\`{e}chet differentiable if this limit is reached uniformly for the pair $x,y$ of points from the unit sphere.

Giles proved in \cite{giles} that
An s.i.p. space is a continuous (uniformly continuous) s.i.p. space if and only if the norm is G\^{a}teaux (uniformly Fr\`{e}chet) differentiable.

From the geometric point of view we know that if $K$ is a $0$-symmetric, bounded,
convex body in the Euclidean $n$-space $\mathbb{R}^n$ (with fixed
origin O), then it defines a norm whose unit ball is $K$ itself (see
\cite{l-g}). Such a space is called \emph{Minkowski space} or normed linear space. Basic
results on such spaces are collected in the surveys
\cite{martini-swanepoel 1}, \cite{martini-swanepoel 2}, and
\cite{martini}. In fact, the norm is a continuous function which is
considered (in geometric terminology, as in \cite{l-g}) as a gauge
function. Combining  this with the result of Lumer and Giles we get
that a normed linear space can be represented as an s.i.p space.
The metric of such a space (called Minkowski metric), i.e. the distance of
two points induced by this norm, is invariant with respect to translations.

Another concept of  Minkowski space was also raised by H. Minkowski and used in Theoretical Physics and Differential Geometry, based on the
concept of indefinite inner product. (See, e.g., \cite{gohberg}.)
The \emph{indefinite inner product} (i.i.p.) on a complex vector space $V$
is a complex function $[x,y]:V\times V\longrightarrow \mathbb{C}$
with the following properties:
\begin{description}

\item[i1]: $[x+y,z]=[x,z]+[y,z]$,
\item[i2]: $[\lambda x,y]=\lambda[x,y]$ \mbox{ for every } $\lambda \in \mathbb{C}$,
\item[i3]: $[x,y]=\overline{[y,x]}$ \mbox{ for every } $x,y\in V$,
\item[i4]: $[x,y]=0$ \mbox{ for every } $y\in V$ then $x=0$.
\end{description}
A vector space $V$ with an i.i.p. is called an \emph{i.i.p. space}.

We recall, that a subspace of an i.i.p. space is positive (non-negative) if all of its nonzero vectors have
positive (non-negative) scalar squares. The classification of subspaces of an i.i.p. space with respect to the positivity property is also an
interesting question. First we pass to the class of subspaces which are peculiar to i.i.p. spaces, and which have no analogous in the spaces with a
definite inner product. A subspace $N$ of $V$ is called \emph{neutral (or isotropic)} if $[v,v]=0$ for all $v\in N$.
In view of the identity
$$
[x,y]=\frac{1}{4}\{[x+y,x+y]+i[x+iy,x+iy]-[x-y,x-y]-i[x-iy,x-iy]\},
$$
a subspace $N$ of an i.i.p. space is neutral if and only if $[u,v]=0$ for all $u,v\in N$. Observe also that a neutral subspace is non-positive and non-negative at the same time, and that it is necessarily degenerate. Therefore it can be proved that an non-negative (resp. non-positive) subspace is the direct sum of a positive (resp. negative) subspace and a neutral subspace. We note that the decomposition of a non-negative subspace $U$ into a
direct sum of a positive and a neutral component is, in general, not unique. However, the dimension of the positive summand is uniquely
determined.

The standard  mathematical model of space-time is a four dimensional i.i.p. space with signature $(+,+,+,-)$, also called Minkowski space in the literature. Thus we have a well known homonyms with the notion of Minkowski space! In our paper we call \emph{Lorentzian space} the Minkowski space defined by an indefinite scalar product with signature $(+,\ldots,+,-)$.

Let {\bf s1}, {\bf s2}, {\bf s3}, {\bf s4}, be the
four defining properties of an s.i.p., and {\bf s5} be the homogeneity
property of the second argument imposed by Giles, respectively.
(As to the names: {\bf s1} is the additivity property of the first argument,
{\bf s2} is the homogeneity property of the first argument, {\bf s3}
means the positivity of the function, {\bf s4} is the
Cauchy-Schwartz inequality.)

On the other hand,  {\bf i1}={\bf s1}, {\bf i2}={\bf s2}, {\bf i3} is the antisymmetry property and {\bf i4} is the nondegeneracy property of the
product. It is easy to see that {\bf s1}, {\bf s2}, {\bf s3}, {\bf s5} imply {\bf i4}, and if $N$ is a positive (negative) subspace of an
i.i.p. space, then {\bf s4} holds on $N$. In the following definition we combine the concepts of s.i.p. and i.i.p..

\begin{defi}
The \emph{semi-indefinite inner product (s.i.i.p.)} on a complex vector space $V$ is a complex function $[x,y]:V\times V\longrightarrow \mathbb{C}$ with the following properties:
\begin{description}

\item[1] $[x+y,z]=[x,z]+[y,z]$ (additivity in the first argument),
\item[2] $[\lambda x,y]=\lambda[x,y]$ \mbox{ for every } $\lambda \in \mathbb{C}$ (homogeneity in the first argument),
\item[3] $[x,\lambda y]=\overline{\lambda}[x,y]$ \mbox{ for every } $\lambda \in \mathbb{C}$ (homogeneity in the second argument),
\item[4] $[x,x]\in \mathbb{R}$ \mbox{ for every } $x\in V$ (the corresponding quadratic form is real-valued),
\item[5] if either $[x,y]=0$ \mbox{ for every } $y\in V$ or $[y,x]=0$ for all $y\in V$, then $x=0$ (nondegeneracy),
\item[6] $|[x,y]|^2\leq [x,x][y,y]$ holds on non-positive and non-negative subspaces of V, respectively. (the Cauchy-Schwartz inequality is valid on positive and negative subspaces, respectively).
\end{description}
A vector space $V$ with a s.i.i.p. is called an s.i.i.p. space.
\end{defi}

We conclude that an s.i.i.p. space is a
homogeneous s.i.p. space if and only if the property {\bf s3} holds, too. An s.i.i.p. space is an i.i.p. space if and only if the s.i.i.p. is an antisymmetric product. In this latter case $[x,x]=\overline{[x,x]}$ implies {\bf 4}, and the function is also Hermitian linear in its second argument. In fact, we have:
$[x,\lambda y+\mu z]=\overline{[\lambda y+\mu z,x]}=\overline{\lambda} \overline{[y,x]}+\overline{\mu} \overline{[ z,x]}=\overline{\lambda}[x,y]+\overline{\mu}[x,z]$. It is clear that both of the classical "Minkowski spaces" can be represented either by an
s.i.p or by an i.i.p., so automatically they can also be represented as an s.i.i.p. space.

It is possible that the s.i.i.p. space $V$ is a direct sum of its two subspaces where one of them is positive and the other one is negative. Then there are two more structures on $V$, an s.i.p. structure (by Lemma 2) and a natural third one, which was called by Minkowskian  structure.
Let $(V,[\cdot,\cdot])$ be an s.i.i.p. space. Let $S,T\leq V$ be positive and negative subspaces, where $T$ is a direct complement of $S$ with
respect to $V$. Define a product on $V$ by the equality
$$
[u,v]^+=[s_1+t_1,s_2+t_2]^+=[s_1,s_2]+[t_1,t_2],
$$ where $s_i\in S$ and $t_i\in T$,
respectively.  Then we say that the pair $(V,[\cdot,\cdot]^+)$ is a \emph{ generalized Minkowski space with Minkowski product $[\cdot,\cdot]^+$}. We also say
that $V$ is a \emph{ real generalized Minkowski space} if it is a real vector space and the s.i.i.p. is a real valued function.

The Minkowski product defined by the above equality satisfies properties {\bf 1}-{\bf 5} of the s.i.i.p.. But in general, property {\bf 6} does not hold. (See an example in \cite{gho 1}.)

If now we consider the theory of s.i.p in the sense of Lumer-Giles, we have a natural concept of orthogonality.  The vector \emph{ $x$ is orthogonal to the vector $y$} if $[x,y]=0$.
Since s.i.p. is neither antisymmetric in the complex case nor symmetric in the real one, this definition of orthogonality is not symmetric in general.

Let $(V,[\cdot,\cdot])$ be an s.i.i.p. space, where $V$ is a complex (real) vector space. The orthogonality of such a space can be defined an analogous way to the definition of the orthogonality of an i.i.p. or s.i.p. space.
The vector \emph{ $v$ is orthogonal to the vector $u$} if $[v,u]=0$. If $U$ is a subspace of $V$, define the orthogonal companion of $U$ in $V$ by
$$
U^\bot =\{v\in V | [v,u]=0 \mbox{ for all } u\in U\}.
$$

We note that, as in the i.i.p. case, the orthogonal companion is always a subspace of $V$. The orthogonal companion of a non-neutral vector $u$ is a subspace having a direct complement of the linear hull of $u$ in $V$. The orthogonal companion of a neutral vector $v$ is a degenerate subspace of dimension $n-1$ containing $v$.

Let $V$ be a generalized Minkowski space. Then we call a vector \emph{ space-like, light-like, or time-like} if its scalar square is positive, zero, or
negative, respectively. Let $\mathcal{S}, \mathcal{L}$ and $\mathcal{T}$ denote the sets of the space-like, light-like, and time-like vectors,
respectively.

In a finite dimensional, real generalized Minkowski space with $\dim T=1$ these sets of vectors can be characterized in a geometric way. Such a space
is called by a \emph{generalized space-time model}. In this case $\mathcal{T}$ is a union of its two parts, namely
$$
 \mathcal{T}=\mathcal{T}^+\cup \mathcal{T}^-,
$$
where
$$
\mathcal{T}^+=\{s+t\in \mathcal{T} | \mbox{ where } t=\lambda e_n \mbox{ for } \lambda \geq 0\} \mbox{ and }
$$
$$\mathcal{T}^-=\{s+t\in \mathcal{T} |
\mbox{ where } t=\lambda e_n \mbox{ for } \lambda \leq 0\}.
$$

It has special interest, the imaginary unit sphere of a finite dimensional, real, generalized space-time model. (See Def.8 in \cite{gho 1}.) It was given a metric on it, and thus got a structure similar to the hyperboloid model of the hyperbolic space embedded in a
space-time model. In the case when the space $S$ is an Euclidean space this hypersurface is a model of the $n$-dimensional hyperbolic space thus it is such-like generalization of it.

In \cite{gho 1} was proved the following statement: $\mathcal{T}$ is an open double cone with boundary $\mathcal{L}$, and the positive part $\mathcal{T}^+$ (resp. negative part $\mathcal{T}^-$) of $\mathcal{T}$ is convex.

We note that if $\dim T> 1$ or the space is complex, then the set of time-like vectors cannot be divided into two convex components. So we have to consider that our space is a generalized space-time model.

\begin{defi}
The set $H:=\{ v\in V | [v,v]^+=-1\}$
is called the \emph{ imaginary unit sphere } of the generalized space-time model.
\end{defi}

With respect to the embedding real normed linear space $(V,[\cdot,\cdot]^-)$ (see Lemma 2) $H$ is a generalized two sheets hyperboloid corresponding to the two pieces of $\mathcal{T}$, respectively. Usually we deal only with one sheet of the hyperboloid, or identify the two sheets projectively. In this case the space-time component $s\in S$ of $v$ determines uniquely the time-like  one, namely $t\in T$. Let $v\in H$ be arbitrary. Let $T_v$ denote the set $v+v^\bot$, where $ v^\bot$ is the orthogonal complement of $v$ with respect to the s.i.i.p., thus a subspace.

The set $T_v$ corresponding to the point $v=s+t\in H$ is a positive, (n-1)-dimensional affine subspace of the generalized Minkowski space
$(V,[\cdot,\cdot]^+)$.
Each of the affine spaces $T_v$ of $H$ can be considered as a semi-metric space, where the semi-metric arises from the Minkowski product restricted
to this positive subspace of $V$. We recall that the Minkowski product does not satisfy the Cauchy-Schwartz inequality. Thus the corresponding distance function does not satisfy the triangle inequality. Such a distance function is called in the literature semi-metric (see \cite{tamassy}). Thus, if the set $H$ is sufficiently smooth, then a metric can be adopted for it, which arises from the restriction of the Minkowski product to the tangent spaces of $H$.

\subsection{Hyperbolic and spherical planes as embedded manifolds}

In this subsection we give a natural connection between the metrics of the spherical and hyperbolic planes. It is based on their common trigonometry.  In fact, as early as 1766 Lambert in \cite{lambert} observed that if there are at least two distinct lines
through a given point that don't intersect a given line, then the area of a triangle with angles a,b,c would be $-R^2 (a+b+c-\pi)$ for
some constant $R$.  He knew also the fact that the area of a triangle on a real sphere of radius $R$ is $R^2 (a+b+c-\pi)$. Comparing the two formulas he noted the following:  "one could almost
conclude that the new geometry would be true on a sphere of imaginary radius".  It turns out that if we substitute the distance $a$ by the imaginary distance $ia$  in any trigonometric formula of the spherical geometry we get a corresponding trigonometric formula valid in hyperbolic geometry. While in the spherical plane all elements of Euclidean trigonometry were reviewed in the nineteenth century (see \cite{casey}) the analogous statements in hyperbolic geometry did not were investigated systematically. Some recent papers try to make up this deficiency (see e.g. \cite{gho 2}, \cite{gho 3}, \cite{gho 4}). We also can find interesting results in this direction in another recent research. I would like to mention those papers of N.J. Wilderberg in which he considered trigonometry in universal geometries. As he said in \cite{wildberger} \emph{"The natural connection with the geometry of Lorentz, Einstein and Minkowski comes from a projective point of view, with trigonometric laws that extend to ‘points at infinity’, here called ‘null points’, and beyond to ‘ideal points’ associated to a hyperboloid of one sheet."}. We note that Wildberger's method \emph{"...works over a general field not of characteristic two, and the main laws can be viewed as deformations of those from planar rational trigonometry..."}. We propose for study the book \cite{wildberger 2} in which first arose the thought of "universal geometry" and the paper \cite{wildberger 3} containing a detailed list of Wildberger's work in this direction.

At the end of the previous paragraph we gave a general approach to the imaginary unit sphere of a generalized space-time model. In the easiest situation it leads to the discover of the analytic connection between the spherical plane and the hyperbolic plane. Let denote by $\langle \cdot,\cdot \rangle$ is the Euclidean inner product of the Euclidean space and denote by $[\cdot,\cdot]$ the indefinite inner product of the Lorentzian space of dimension $3$. We can compare the distances of the spherical and hyperbolic geometry because the first distance realizes as a distance of points of the sphere of the Euclidean space of radius $r$, and the latter one as a distance of points in the imaginary sphere of radius $ir$ of the Lorentzian space, respectively. The angle measure of the two rays ($PX$ and $PY$) in these planes is nothing else as the dihedral angle measure between the planes with unit normals $\zeta$ and $\xi$ through the origin intersecting the required sphere in the given rays. From this we get the following results:
$$
\begin{array}{|l|c|c|}
\hline
    & S^2(\mathbb{R}), R=r & H^2(\mathbb{R}), R=ir  \\
\hline
    \cos \frac{\rho(X,Y)}{R} & \frac{\langle x,y\rangle}{R^2} & \frac{[x,y]}{R^2} \\
\hline
   \cos (XPY)\angle & \langle \zeta,\xi \rangle & -[\zeta,\xi]\\
\hline
\end{array}
$$

\subsection{Diagrams on the connections}

\begin{figure}[ht]
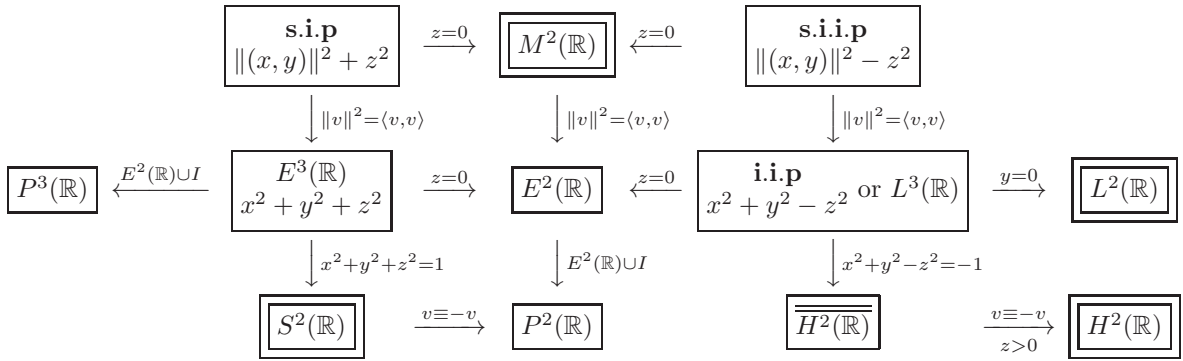

$$
\minCDarrowwidth20pt\begin{CD}
@.\fbox{$\begin{matrix}{\bf s.i.p} \\ \|(x,y)\|^2+z^2 \end{matrix}$} @>z=0>> \fbox{\fbox{$M^2(\mathbb{R})$}} @<z=0<< \fbox{$\begin{matrix} {\bf s.i.i.p}\\ \|(x,y)\|^2-z^2 \end{matrix}$}  \\
@. @VV\|v\|^2=\langle v,v\rangle V @VV\|v\|^2=\langle v,v\rangle V  @VV\|v\|^2=\langle v,v\rangle V   \\
\fbox{$P^3(\mathbb{R})$}@<E^2(\mathbb{R})\cup I<<\fbox{$\begin{matrix}E^3(\mathbb{R})\\ x^2+y^2+z^2 \end{matrix}$} @>z=0>> \fbox{$E^2(\mathbb{R})$} @<z=0<< \fbox{$\begin{matrix} {\bf i.i.p}\\ x^2+y^2-z^2 \end{matrix}$ or $L^3(\mathbb{R})$} @>y=0>> \fbox{\fbox{$L^2(\mathbb{R})$}}  \\
@. @VVx^2+y^2+z^2=1V   @VVE^2(\mathbb{R})\cup IV @VVx^2+y^2-z^2=-1V \\
@. \fbox{\fbox{$S^2(\mathbb{R})$}} @>v\equiv -v>>\fbox{$P^2(\mathbb{R})$}  @. \fbox{$\overline{\overline{H^2(\mathbb{R})}}$} @>v\equiv -v>{z>0}> \fbox{\fbox{$H^2(\mathbb{R})$}}
\end{CD}
$$
\caption{Non-Euclidean planes with associated vector space}
\label{diagram1}
\end{figure}

Our first diagram (Fig. \ref{diagram1}) shows that our general linear algebraic terminology how leads to the four non-Euclidean planes. Using the notation of the first subsection of the introduction {\bf s.i.i.p} is an abbreviation of a semi-indefinite inner product space which is a generalized space-time model of dimension $3$. We also use a coordinate frame of the embedding space, the coordinates can be considered with respect to this system. The used product either an indefinite inner (scalar) product or a semi inner product with the respective abbreviations {\bf i.i.p} or {\bf s.i.p}. An equality over an arrow means the restriction of the space to the set of those vectors which are satisfying it, with the same product. For example the {\bf s.i.i.p} restricted by the equality $z=0$ will be an {\bf s.i.p} showing that the getting plane is a Minkowski (normed) plane.

On the second figure (Fig. \ref{diagram2}) we can see the so-called Cayley-Klein geometries. Three geometries (spherical, hyperbolic and Lorentzian) from the investigated ones belong to this class of geometries. If we consider the real projective plane and fix a group of isometries as a special subgroup of the group of projective collinearities then we can associate a projective metric to this subgroup. This metric determines a plane geometry. By Klein's program these geometries classified in an algebraic way. The classes depend on the type of the measure of the length and also on the type of the measure of the angle. There are three possibilities for both of the measures, they could be elliptic, parabolic or hyperbolic, respectively. On this way we can get nine possibilities which can be realized, these are the so-called Cayley-Klein geometries. The Cayley-Klein geometries have a synthetic description (see in \cite{struve}) and also a linear algebraic one (see the paper of Wilderberg \cite{wildberger} and a recent paper of Juh\'asz \cite{juhasz}).
\begin{figure}[ht]
\renewcommand{\arraystretch}{1.5}
\begin{tabular}{|c|c|c|c|}
  \hline
  angle $\setminus$ length & elliptic & parabolic & hyperbolic \\ \hline
  elliptic & elliptic plane & $E^2(\mathbb{R})$ & $H^2(\mathbb{R})$  \\ \hline
  parabolic & $\overline{E^2(\mathbb{R})}$ & $G^2(\mathbb{R})$ & $\overline{L^2(\mathbb{R})}$ \\ \hline
  hyperbolic & $\overline{H^2(\mathbb{R})}$ & $L^2(\mathbb{R})$ & $\overline{\overline{H^2(\mathbb{R})}}$\\
  \hline
\end{tabular}
\caption{Cayley-Klein geometries}
\label{diagram2}
\end{figure}

The upper line means the projective dual of that geometry which lies under the line. Struve gave a synthetic axiomatic representations for eight geometries from the above nine. For this purpose he defined the dual of the known three axioms on parallels, (so that every two lines are intersecting, holds the Euclidean axiom of parallels or hold the hyperbolic axiom of parallels, respectively) and prove that those are enough to determine the geometries in question. The only geometry which cannot determine in that way is the double hyperbolic space $\overline{\overline{H^2(\mathbb{R})}}$ which is in the associated i.i.p space is the set of those points whose scalar squares are equal to $-1$ with the corresponding metric. This plane is called also by anti-de Sitter space of type $(2,1)$ containing the two branches of the getting hyperboloid. If we identifies the two branches or consider only one of them we get a model for the hyperbolic plane. It is in a strong analogy with the case of the spherical--elliptic pair of planes.

The Cayley-Klein geometries which cannot find in the first diagram are in Fig. \ref{diagram3}. This diagram shows that these also can be get from the embedding indefinite inner product space and an extraction of the diagram in Fig. \ref{diagram1} can contains them, too.

\begin{figure}[ht]
$$
\minCDarrowwidth20pt\begin{CD}
@.\fbox{$\overline{L^2(\mathbb{R})}$}\\
@. @Ax^2+y^2-z^2=-1Az\neq 1A \\
\fbox{$\overline{H^2(\mathbb{R})}$} @<x^2+y^2-z^2=1<< \fbox{$\begin{matrix} {\bf i.i.p}\\ x^2+y^2-z^2 \end{matrix}$} @>y=z>\rho(v_1,v_2):=|y_2-y_1|> \fbox{$G^2(\mathbb{R})$} \\
\end{CD}
$$
\caption{}
\label{diagram3}
\end{figure}

Using Struve approach we now give a short description of those geometries which are not well-known for a reader unfamiliar in non-Euclidean geometries.

The so-called \emph{co-hyperbolic plane} $\overline{H^2(\mathbb{R})}$ means a geometry in which every two lines have a common point of intersection and if $A$ is a point and $a$ is a line not incident with $A$ then there are precisely two points $X$ and $Y$ on $a$ which is parallel to $A$ in the hyperbolic meaning. Two lines are called h-parallel if they have no common point and no common perpendicular. Dually, two points are called h-parallel if they have no joining line and no common polar point. We get a model of this plane if we take the complement of a Cayley-Klein model of the hyperbolic plane with respect to the embedding projective plane. We call here line a projective line which avoid the given hyperbolic model. It can be seen easily that if from the point $A$ we draw tangents to the Cayley-Klein model and consider the intersections of these tangents with the line $a$ get the points $X$, $Y$ are parallels to $A$. There are further two names of this plane it can be called as \emph{"hyperbolic plane with positive curvature"} (see \cite{romakina}), or also \emph{de Sitter space of dimension $2$} (see in \cite{gho 7}). With indefinite inner product we can model it with a one-sheet hyperboloid,  containing those points of the space whose coordinates fulfill the equality $x^2+y^2-z^2=1$.

The incidence structure of a \emph{co-Euclidean plane} $\overline{E^2(\mathbb{R})}$ can be modelled by the removal of a pencil of lines (with its base) from the  projective plane $\overline{P^2(\mathbb{R})}$. Hence a metric model for it cannot be get from the indefinite inner product space of dimension 3 even from the Euclidean space of dimension 3 (which is also a semi inner product space). We can consider the elliptic geometry modeled by the unit sphere from which we remove one of its points and redefined the set of lines to the set of those lines which do not through the points removed from the model. We can also consider the metric of the elliptic plane to this restricted sets of points and lines.

$G^2(\mathbb{R})$ is called by \emph{Galilean plane}. $G^2(\mathbb{R})$ is an affine space so it can be embedded into the i.i.p space of dimension three. To get a model for this geometry we have to remove a line and also a pencil of rays from the embedding projective plane. Hence it can be demonstrated as the geometry of a suitable Euclidean plane of the i.i.p space containing precisely one vectors with zero length. An appropriate choice to get a model if we consider the bisector of the coordinate planes $(x,y)$ and $(x,z)$ with that degenerated metric which is based on the function $\rho(v_1,v_2)=|y_1-y_2|$ to measure the distance of points (see \cite{artikbayev}).

The incidence structure of the \emph{co-Lorentzian plane} $\overline{L^2(\mathbb{R})}$ can be get by the removal of one pencil of rays from the hyperbolic plane. A natural model of it the hyperboloid model of the hyperbolic plane without its intersection point $(0,0,1)$ with the $z$ axis. The line set contains those hyperbolic lines which are not goes through this point and the metric is the metric of the hyperbolic plane. This plane is called by \emph{quasi-hyperbolic plane} in isotropic geometry (see \cite{sacks}, \cite{yaglom}, \cite{slipcevic}).

\section{Conics}

\subsection{Hyperbolic and spherical conics}
There are several papers on conics in hyperbolic or spherical planes, respectively. Probably, the earliest complete list with respect to the hyperbolic plane can be found in a work of the Hungarian scientist Cyrill V\"or\"os who wrote a nice book on analytic hyperbolic geometry in Hungarian \cite{voros}.

In the second half of the previous century  Emil Moln\'ar gave a nice classification with a synthetic approach (see in \cite{molnar}). This approach based on reflections at lines. (Note that these isometries generate the full group of isometry). A conic $\gamma$ in the sense of Moln\'ar, with foci
$A,B\in P^2(\mathbb{R})$, is defined by choosing a line $x_1$ in $P^2(\mathbb{R})$ which is not a boundary line (i.e., $x_1$ is not tangent to the absolute) and not passing through either $A$ or $B$, and with $A$ and $B$ not each other's reflections across $x_1$. Then $\gamma$ consists of
points $X_{11}$ and $X$ chosen as follows. $X_{11}$ is the intersection of the lines $a_{11}$ through
$A$ and the reflection $Bx_1$ of $B$ across $x_1$ and $b_{11}$ through $B$ and $Ax_1$ . (The line $x_1$ is chosen so that $a_{11}$ and $b_{11}$ are not boundary lines.) The other points $X$ are defined by fixing a point $Y$ on $x_1$ and taking the lines $a$ through $Y$ and $A$ and $b$ through $Y$ and $B$, and then if neither $a$ nor $b$ is a boundary line, letting $X$ be the intersection of $a^a_{11}$ and $b^b_{11}$ (the reflections of $a_{11}$ and $b_{11}$ across $a$ and $b$, respectively). This definition based on that property of Euclidean conics that the line determined by the point of intersection of two tangents with a focus, is the bisector of the focal radii corresponding to those points of the conic which are on the given tangents, respectively. (It can be seen easily that this metric property has a projective interpretation, too.) This definition leads to the same types of conics which can be get using the fact that with respect to the Cayley-Klein model every hyperbolic conic is the intersection curve of the model disk with a projective conic of the embedding projective plane.

We can find also two papers of K.Fladt (\cite{fladt1},\cite{fladt2}) containing a complete analytic classification. This latter work inspired a characterization with dual pairs of conics by G. Csima and J. Szirmai in \cite{csima}. Interesting problem that what does it means the phrase "conic section"? Chao and Rosenberg (\cite{chao}) wrote a paper on the hyperbolic concepts of conics giving the logical equivalence and non-equivalences among them.

We also have to mention (without the claim of completeness) some further references on hyperbolic conics. The paper of G. Weiss \cite{weiss} contains interesting metric definitions and theorems working in the hyperbolic and elliptic planes, respectively. There are two papers by H.P. Schr\"ocker and M.J. Weber on the minimal area ellipses in elliptic and hyperbolic plane, respectively (see \cite{schroecker-weber 1}, \cite{schroecker-weber 2}).

On spherical conics we can find the earliest paper of Sykes and Pierces \cite{sykes} at the end of the nineteenth century. We mention here some other papers with similar results written by Dirnb\"ock \cite{dirnbock} at the end of the last century the recent paper of Altunkaya at all. from 2014 \cite{altunkaya} and the nice book of Glaeser, Stachel and Odehnal \cite{glaeser,stachel}. This latter contains valuable informations of non-Euclideans conics, too. Chapter 10.1 contains a nice approach to spherical conics with several results explained with respect to spherical geometry and compared with the Euclidean situations.

\subsubsection{Classification of spherical conics}

We use here the approach of the paper of Sykes and Pierces.
A \emph{spherical conic} is the intersection of a unit-sphere with a cone of the second degree, whose vertex is at the centre of the sphere. Since the cone is double, it will cut the sphere in two closed curves; and we therefore name the conic differently according to the hemisphere considered. If the sphere be divided by the principal plane of the cone, it gives a closed curve whose centre will be the pole of the dividing circle, and whose principal diameters will be the arcs of the greatest and least sections of the cone. This form of conic is a \emph{Spherical Ellipse}. If the sphere be divided by the plane of least section of the cone, the conic will consist of two branches. Its centre will be the pole of the dividing circle, and its principal diameters will be the arcs made by the plane of greatest section of the cone and the principal plane. This curve is the \emph{Spherical Hyperbola}. If, again, the sphere be bisected by a plane perpendicular to the two already mentioned, there is still a third form of spherical conic, having its centre at the pole of the bisecting circle. There is, properly speaking, as might be expected from the method of projection used, no spherical parabola. If a plane parabola be projected upon a sphere, points at infinity are projected, and the spherical parabola is merely an ellipse or an hyperbola. The conic of which the major axis is a quadrant has, however, the closest analogy to the Parabola. We note that the above description of conics (quoted from \cite{sykes}) shows that in spherical geometry there is only one type of conics which can be get also in an analytic manner.

A spherical conic may also be defined as the locus of an equation of the second degree in spherical co-ordinates. The general equation is
$$
ax^2+2hxy +by^2+2gx+2fy+c = 0.
$$
This can be transformed to the centre as origin; and, if we choose the principal diameters as axes, it can be reduced to the form
$$
\frac{x^2}{a^2}+\frac{y^2}{b^2}=1.
$$
The equation for determining the centre is a cubic, and this shows that a spherical conic has three centres.

We refer here the Theorem 10.1.3 from \cite{glaeser} showing the above facts.
\begin{theorem}[\cite{glaeser}]
Any spherical ellipse with focal points $S_1$ and $S_2$ and major axis $2a$ is a branch of a spherical hyperbola with ordered foci $S_2$ and $S_1^{\star}$ and major axis $\pi-2a$, where $S_1^{\star}$ is the antipode of $S_1$.
\end{theorem}
This nice book contains also some interesting theorems working also in non-Euclidean planes as the so-called Ivory's theorem.

\subsubsection{Classification of hyperbolic conics}

We review in this paragraph the classification of conics on the base of their analytic definition. Our originated is the work \cite{csima}. The classification of the conics on the extended hyperbolic plane can be obtained in dual pairs. (On the projective extension of the hyperbolic plane I propose the study of the paper \cite{gho 2}.)

Consider a one parameter conic family of our point conic with the absolute conic, defined by
$$
x^T({\bf a}+\rho {\bf e})x=0.
\notag
$$

Since the characteristic polynomial $\Delta(\rho):=\det({\bf a}+\rho {\bf e})$ is an odd degree one, this conic pencil has at least one real degenerate element $(\rho_1)$, which consists of at most two point sequences with holding lines $\Bsp_1^1$ and $\Bsp_1^2$ called asymptotes. Therefore we get a product
$$
x^T({\bf a}+\rho_1 {\bf e})x=(\Bsp_1^1x)^T(\Bsp_1^2x)=x^T((\Bsp_1^1)^T\Bsp_1^2)x=0
\notag
$$
with occasional complex coordinates of the asymptotes. Each of these two asymptotes has at most two common points with the absolute and with the conic as well. Thus, the at most 4 common points with at most 3 pairs of asymptotes can be determined through complex coordinates and elements according to the at most 3 different eigenvalues $\rho_1$, $\rho_2$ and $\rho_3$.

In complete analogy with the previous discussion in dual formulation we get that the one parameter conic family of a line conic with the absolute has at least one degenerate element $(\sigma^1)$ which contains two line pencils at most with occasionally complex holding points $\bmf_1^1$ and $\bmf_2^1$ called foci.
$$
u({\bf A}+\sigma^1 {\bf E})u^T=(u \bmf_1^1)(u \bmf_2^1)^T=u(\bmf_1^1(\bmf_2^1)^T)u^T=0
\notag
$$
For each focus at most two common tangent line can be drawn to the absolute and to the line conic. Therefore, at most four common tangent lines with at most three pairs of foci can be determined maybe with complex coordinates to the corresponding eigenvalues $\sigma^1$, $\sigma^2$ and $\sigma^3$.

Combining this discussions with in \cite{fladt1} the classification of the conics on the extended hyperbolic plane can be obtained in dual pairs.

To find an appropriate transformation, so that the resulted normalform characterizes the conic, we take a rotation around the origin $O(0,0,1)^T$ and a translation parallel to $x^2=0$. The characteristic equation
$$
\Delta(\rho)=\det({\bf a}+\rho{\bf e})=
\det\left(\begin{array}{ccc}
	a_{11}+\rho & a_{12} & a_{13}\\
	a_{21} & a_{22}+\rho & a_{23}\\
	a_{31} & a_{32} & a_{33}-\rho
\end{array}\right)=0
$$
has at least one real root denoted by $\rho_1$.

This is helpful to determine the exact transformation if  the equalities $\rho_1=\rho_2=\rho_3$  not hold.
With this transformations we obtain the normalform
$$
\rho_1 x^1x^1+a_{22}x^2x^2+2a_{23}x^2x^3+a_{33}x^3x^3=0.
\label{nf}
$$
In the following we distinguish $3$ different cases according to the other two roots:

\begin{itemize}
\item \emph{The case of two different real roots:}
Then the monom $x^2x^3$ can be eliminated from the equation above, by translating the conic parallel to $x^1=0$. The final form of the conic equation in this case, called \emph{central conic} section:

\[\rho_1 x^1x^1+\rho_2 x^2x^2-\rho_3 x^3x^3=0.\]

Because the conic is non-degenerate $\rho_3x^3\neq0$ follows and with the notations $a=\frac{\rho_1}{\rho_3}$ and $b=\frac{\rho_2}{\rho_3}$ the matrix can be transformed into ${\bf a}=diag\left\{a,b,-1\right\}$, where $a\leq b$ can be assumed. The equation of the dual conic can be obtained using the polarity ${\bf E}$ respected to the absolute by ${\bf E}~{\bf A}~{\bf E}^{-1}=diag\{\frac{1}{a},\frac{1}{b},-1\}$.
By the above considerations we can give an overview of the generalized central conics with representants:

If the conic section has the normalform $a x^2+b y^2=1$ then we get the following types of central conic sections:
\begin{center}
\begin{tabular}{|l|r|}
\hline
	Absolute conic & $a=b=1$ \\
\hline
	Circle & $1<a=b$ \\
    Circle enclosing the absolute &$a=b<1$ \\
\hline
	Hypercycle & $1=a<b$ \\
    Hypercycle enclosing the absolute & $0<a<1=b$  \\
\hline
	Hypercycle excluding the absolute & $a<0<1=b$ \\
	Concave hyperbola &  $0<a<1<b$ \\
\hline
	Convex hyperbola &   $a<0<1<b$ \\
	Hyperbola excluding the absolute & $a<0<b<1$ \\
\hline
	Ellipse & $1<a<b$ \\
	Ellipse enclosing the absolute &  $0<a<b<1$\\
\hline
	empty& $a\leq b\leq0$ \\
\hline						
\end{tabular}
\end{center}
where the conic and its dual pair lies in the same row. (If it is not self-dual we have to take the parameter transformations $a'=\frac{1}{a}$ and $b'=\frac{1}{b}$).

\item \emph{The case of coinciding real roots:}
The last translation cannot be enforced but it can be proved that $\rho_2=\rho_3=\frac{a_{33}-a_{22}}{2}$ follows.
With some simplifications of the formulas in \cite{fladt1} we obtain the normalform of the so-called generalized \emph{parabolas}.
The parabolas have the normalform $ax^2+(b+1)y^2-2y=b-1$ and the following cases arise:
\begin{center}
\begin{tabular}{|l|r|}
\hline
				Horocycle & $0<a=b$ \\
				Horocycle enclosing the absolute &  $a=b<0$ \\
    \hline
				Elliptic parabola & $0<b<a$ \\
				Parabola enclosing the absolute & $b<a<0$ \\
	\hline
			    Two sided parabola & $a<b<0$ \\
				Concave hyperbolic parabola &  $0<a<b$ \\
	\hline
				Convex hyperbolic parabola & $a<0<b$ \\
				Parabola excluding the absolute  & $b<0<a$ \\
	\hline
\end{tabular}
\end{center}
where the dual pairs with parameters $a'=-\frac{b^2}{a}$ and $b'=-b$ are the same row on the table.

\item \emph{The case of two conjugate complex roots:}\\
Then the last translation cannot be performed to eliminate the monom $x^2x^3$ but we can eliminate the monom $x^3x^3$ by an appropriate transformation described in \cite{fladt1}. Shifting to inhomogeneous coordinates and simplifying the coefficients we obtain:
The so-called \emph{semi-hyperbola} has the normalform $a x^2+2 b y^2-2y=0$ where $\left|b\right|<1$ and its dual pair is projectively equivalent with another semi-hyperbola with $a'=\frac{1}{a}$ and $b'=-b$.

\item Overviewing the above cases only one remains, when the conic has no symmetry axis at all and $\rho_1=\rho_2=\rho_3$.
Ignoring further explanations we claim the following:
If the conic has the normalform
$(1-x^2-y^2)+2a y(x+1)=0$ where $a>0$ then it is called \emph{osculating parabola}. Its dual is also an osculating parabola by a convenient reflection.
\end{itemize}

The following theorem summarize the above information:

\begin{theorem}{\cite{csima}}
There are $22$ types of hyperbolic conics listed in the above discussion. Most of them are self-dual and the other ones can be paired into dual-pairs.
\end{theorem}

\subsection{Conics in Lorentzian-plane}
In a Lorentzian plane also there are two possibilities (analytic and metric ones, respectively) to define a conic. In the paper of Birkhoff and Morris  \cite{birkhoff} the definition based on the two foci property and hence each of the types of conics have a singular definition. In such a way we can distinguish curves not only on their metric properties but on the base of their causal characters, e.g. we can say about relativistic time-like ellipse or relativistic time-like hyperbola, respectively. More precisely, if $T(E,E')= \sqrt{(t-t')^2-(x-x')^2}$ means the so-called time interval between two events $E(t,x)$ and $E'(t',x')$ then for fixed positive $a$ we consider the relativistic \emph{time-like ellipse}, as the locus of all points $E(t,x)$ which satisfy
$$
T(F,E)+T(F',E)=2a
$$
with respect to two fixed point $F$ and $F'$ so that their segment is time-like so that $T(F,F')=2c$ is real. These points are the foci of the ellipse. Similarly defined the two branch of a \emph{relativistic hyperbola} by the equalities
$$
T(F,E)-T(F',E)=2a \mbox{ and } T(F',E)-T(F,E)=2a,
$$
respectively. It has been proved, that each geometrical conic of the form
$$
\frac{t^2}{a^2}-\frac{x^2}{a^2-c^2}=1 \quad a>0,c>0,a\ne c
$$
is the union of pieces consisting of confocal relativistic conics. (For more details see the paper \cite{birkhoff}.) It proved also that the relativistic conics are geometrically tangent to the null (isotropic) lines through the foci $F$ and $F'$.

Similar definitions can be found in the Shonoda's paper \cite{shonoda} (and also in the recent arxive \cite{sanchez}). In this paper the Apollonius definition of conics were used to generate algebraic curves in the Minkowski space-time plane (in the Lorentzian plane). It turn out to be different from classical conic sections. It has been extended and classified the sort of “M-conics”. Also discussed the cases of the singularity points of these M-conics, coming from the transition from time-like world to space-like world through the light-like one.

In \cite{stachel}, \cite{stachel1} and \cite{gho 5} the authors used linear algebra to the definition. In fact, from a geometric point of view, the conics can be represented by a quadric as the zero set of a quadratic form (or the zero set of its symmetric bilinear form). Fixing a regular symmetric bilinear form as an indefinite inner product $\langle \cdot,\cdot \rangle$ any quadric can be regarded as the zero set of a symmetric bilinear function $\langle x, l(y)\rangle$, where $l$ is a selfadjoint transformation with respect to the product. Since in an inner product space the  self-adjoint transformations has a classification (see \cite{gohberg}) it has also a classification of conics defined by in this way.

There is a possibility to define conics in Lorentzian plane by projective geometry, since it is also an affine Cayley-Klein geometry. The Lorentzian plane is a projective plane where the metric is induced by a real line $f$ and two real points $F_1$ and $F_2$ incidental with it. A curve in the Lorentzian plane is circular if it passes through at least one of the absolute points. If it does not share any point with the absolute line except the absolute points, it is said to be entirely circular. Every curve of order $n$ intersects the absolute line in $n$ points. If $F_1$ is an intersection point of the curve and the absolute line with the intersection multiplicity $r$ and $F_2$ is an intersection point of the curve and the absolute line with the intersection multiplicity $t$, then it is said to be a curve with the type of circularity $(r, t)$ and its degree of circularity is defined as $r + t$. The classification with respect to this method can be found in the papers \cite{kovacevic} and \cite{jurkin} and says:
\begin{theorem}[\cite{kovacevic}]
The conics are classified into: non-circular conics (ellipses, hyperbolas, parabolas), special hyperbolas (circularity of type $(1, 0)$), special parabolas (circularity of type $(2, 0)$) and circles (circularity of type $(1, 1)$).
\end{theorem}

The subsection 10.2 in \cite{glaeser,stachel} contains also a nice description of the above six types of conics. In addition, the last pages of the book contains also the definition of hyperbolic conics as conics which are projective conics of the embedding projective plane containing interior points of the Cayley-Klein model. This definition needed to state that the Apollonian definition of conics (based on the property of the leading line and a focus) does not work in hyperbolic plane. (In my point of view this metrical property gives another types of curves in the hyperbolic plane as the analytic one.) Further results on conics in Lorentzian plane can be found in \cite{wildberger 4}.

\subsection{Conics in Minkowski normed planes}

While the projective geometry under the Minkowski plane is the same as under the Euclidean plane (it is an affine plane, too) the used metric kills the possibility that we classify projective conics using metric properties. From this reason we have several immediate definitions for conics. The first papers in this direction is the paper of Wu, Ji and Alonso \cite{wu}, and the paper of G.Horv\'ath and Martini  \cite{gho 5}, respectively. Minkowski conics was investigated by also Fankhanel in the paper \cite{fankhanel}. In normed planes we have three different possibilities to define ellipses metrically. The first one was investigated in the paper \cite{wu}. In \cite{gho 5} we can find the following definitions refer to a normed plane $X$.

\begin{itemize}
\item
Let ${\bf x},{\bf y}\in X$, ${\bf x}\not ={\bf y}$, and $2a\geq 2c=\|{\bf x}-{\bf y}\|$.  The set
$$
E({\bf x},{\bf y}, a)=\{{\bf z} \in X: ||{\bf z}-{\bf x}|| +||{\bf z}-{\bf y}||=2a \}
$$
is called \emph{the ellipse defined by its foci} $x$ and $y$.

\item
Let $L:=(2a)\cdot K$ be a homothetic copy of the unit disk $K$, and ${\bf x}\in L$ be an arbitrary point from it. The locus of points ${\bf z}\in X$ for which there is a positive $\varepsilon$ such that ${\bf z}+\varepsilon K$ touches $L$ and contains ${\bf x }$ on its boundary is called \emph{the ellipse defined by its leading circle and its focus} {\bf x}.

\item
Let $l$ be a straight line, ${\bf x}$ a point, and $\gamma =\frac{a}{c}$ a ratio larger than 1. The locus of points ${\bf z}\in X$, for which there is a positive ${\varepsilon}$ such that the boundary of the disk ${\bf z}+\varepsilon K$ contains {\bf x} and the disk ${\bf z}+\gamma (\varepsilon K)$ touches the line $l$, is called \emph{the ellipse defined by its leading line and its focus} {\bf x}.
\end{itemize}

\begin{theorem}[\cite{gho 5}]
In any normed plane the following holds: an ellipse, defined by its foci, is always an ellipse defined by its leading circle and a focus, and the converse statement is also true. On the other hand, an ellipse defined by its leading line and a focus is not necessarily an ellipse defined by its foci, and again the converse is true.
\end{theorem}

\begin{figure}[ht]
 \centering
  \includegraphics[height=5cm]{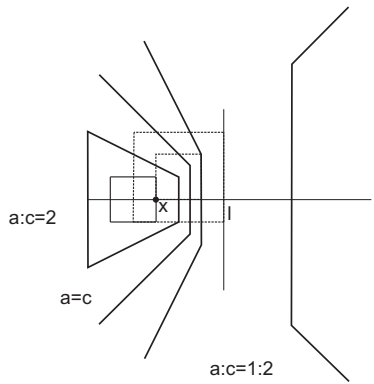}\\
  \caption{Conics on the $l_\infty$ plane}
  \label{mc30}
\end{figure}

In Fig.\ref{mc30} we can see
that there is an ellipse following the third definition which is not centrally symmetric. By Theorem 2 of \cite{wu} it is not an ellipse by the first definition. In our example the norm is the $L_\infty$ norm, and the leading line $l$ and the focus ${\bf x}$ are in ``symmetric position'' with respect to the circle of this Minkowski plane, which is a square.

\begin{figure}[ht]
    \centering
  \includegraphics[height=5cm]{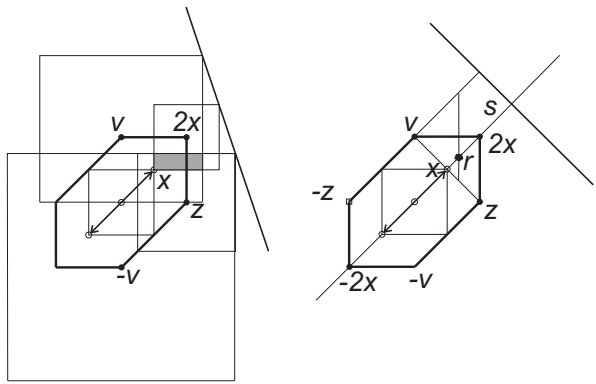}\\
  \caption{A metric ellipse which has no leading line}
  \label{mc40}
\end{figure}

Conversely, consider the ellipse $E(-{\bf x},{\bf x}, 2)$ defined by its foci and shown in Fig.\ref{mc40}.
First we can see that if it is also an ellipse defined by its leading line, then the leading line $l$ and the new focus ${\bf x}'$ have to be in ``symmetric position'' with respect to the line joining the original foci. ``Symmetric'' means that this line is parallel to a diagonal of the unit square. In fact, if this is not the case, we get a figure as shown on the left side of Fig.\ref{mc40}. The squares $S_{2{\bf x}}, S_{\bf v}, S_{\bf z}, S_{-{\bf v}}$ with centers $2{\bf x}$, ${\bf v}$, ${\bf z}$, $-{\bf v}$, respectively, touch $l$. The focus has to lie in the shaded rectangle, as the common point of the boundaries of homothetic copies $2{\bf x}+\frac{c}{a}S_{2{\bf x}}$, ${\bf v}+\frac{c}{a}S_{{\bf v}}$ and ${\bf z}+\frac{c}{a}S_{{\bf z}}$ of such squares (with a homothety ratio smaller than 1). On the other hand, the boundary of the square $-{\bf v}+\frac{c}{a}S_{-{\bf v}}$ intersects the shaded rectangle in a segment parallel to that one in which it is intersected by ${\bf z}+\frac{c}{a}S_{{\bf z}}$. So it is impossible to give a
good position for the focus $x′$.

We now assume that $l$ and ${\bf x}'$ have symmetric position (see the right side of Fig.\ref{mc40}).
If this holds and the Euclidean distance of $l$ and $2{\bf x}$ is $s$, and that of ${\bf x'}$ and ${\bf x}$ is $r$, then, using the fact that the points $2{\bf x}$, $-2{\bf x}$ and ${\bf v}$
have to lie on the new ellipse, we have the equalities
$$
\frac{r}{s}=\frac{4-r}{4+s}=\frac{2-r}{1+s}\,,
$$
implying that
$$
s=1 \mbox{ and } r=\frac{2}{3}
$$
and showing that $\frac{a}{c}=\frac{2}{3}$.
Thus the leading line and the focus are both determined. On the other hand, the point $-{\bf z}$ is not on the obtained ellipse, since the required ratio for it is $\frac{12-\sqrt{2}}{12}\not =\frac {2}{3}$.

For hyperbolas there are similar metric definitions. It can be proved that in normed planes, a hyperbola defined by its foci is always a hyperbola defined by its leading circle and a focus. The converse statement is also true. In general, the third definition yields a different class of curves.

For the case of parabolas, the first two definitions have no analogue, and so we have only the third case.
In a normed plane, let $l$ be a straight line, and ${\bf x}$ be a point. The locus of the points ${\bf z}\in S$ for which there is a positive ${\varepsilon}$ such that the boundary of the disk ${\bf z}+\varepsilon K$ contains {\bf x} and touches the line $l$, will be called \emph{the parabola defined by its leading line and its focus} {\bf x}.

It is also true that the metric parabola is a simple curve which does not contain segments if and only if the normed plane under consideration is strictly convex.

Finally some words about the analytic building up of the projective conics of a normed space. The following way is a possibility to define quadrics in the projective augmentation of any smooth, strictly convex space. We describe this method in the two-dimensional case, where the quadric is clearly a conic.

Every normed plane can be represented as a semi-inner product space (s.i.p.; see \cite{lumer} and \cite{giles}). If the unit disk  is strictly convex, this representation is unique. As proved in \cite{giles}, the orthogonality with respect to the s.i.p. is equivalent to the orthogonality concept of Birkhoff (see, e.g., \cite{alonso1} and \cite{alonso2}).  Koehler proved in \cite{koehler} that if the generalized Riesz-Fischer
representation theorem is valid in a normed space, then every bounded linear operator $A$ has a generalized adjoint $A^T$ defined by the equality
$$
[A(x),y] = [x,A^T(y)] \mbox{ for all } x,y \in V.
$$
It can be proved that if in all strictly convex and smooth spaces the above assumption holds, then in such a space there is a generalized adjoint.
We remark that $A^T$ is in general not a linear transformation. We say that the linear mapping is \emph{self-adjoint} if $A=A^T$.
If $A$ is self-adjoint, then any element of its class in the Projective General Linear Group of $V$ is self-adjoint, too. So we can call such a family of
operators \emph{class of self-adjoint linear operators of the projective space} $P(V)$. Now the concept of conics can be introduced as follows.
Let $P(V)$ be a real projective space with the two-dimensional semi-inner product space $(V,[\cdot,\cdot])$. A (non-degenerate) projective conic is the zero set of a (non-degenerate) form $\Phi(x, y) = [A(x),y]$, with an invertible self-adjoint operator $A$ of $P(V)$.

We remark that the form $\Phi(x,y)$ is linear in its first argument, homogeneous in its second one, but is neither symmetric, bilinear nor positive. It is symmetric and bilinear if the semi-inner product is symmetric; bilinear if the semi-inner product is additive in its second argument; and positive if $A$ is a square operator (meaning that it is the square of another self-adjoint operator, denoted by $\sqrt{A}$).

The group of self-adjoint operators is basically determined by the unit disks, and it determines the projective conics in analytic sense. Thus, in this
setting the metric of the plane is also used for smooth, strictly convex normed planes. We finish with two problems:
\begin{itemize}
\item Characterize the self-adjoint operators for smooth, strictly convex normed planes.
\item Describe relations between metric conics and general ones.
\end{itemize}

The first question was also investigated in \cite{langi} in the case when the plane also has a Lipshitz-type property. Some further observations can be found in \cite{gho 9}.

\section{roulettes}

The investigation of roulettes in geometry is important not only with respect to their nice geometric properties but also their influence on kinematics. Hence first of all we give a short overview of this connection. After this review we consider the validity of the Euclidean results in our non-Euclidean geometries. Similarly to the case of conics in the Minkowski plane we have to define new apparat to consider roulettes while in the cases of the other three geometries the roulettes have an analogous approach as in the Euclidean situation.

\subsection{Motions of rigid systems in the Euclidean plane}

Consider a plane $\Sigma'$ which is moving on the fixed plane $\Sigma $. The two simplest possibilities for such movements are given by translation and rotation. In Euclidean geometry we can substitute the planes with cartesian coordinate frames $Oxy$ and $O'uv$, and when we would like to describe the motion of a point $P$ of the moving plane, we need the coordinates $u,v$ of the point $P$ in the moving frame, the coordinates $p,q$ of $O'$ in the fixed coordinate system, and the angle $\varphi $ of the positive half of the $X$-axis of the fixed frame with the positive half of the $x$-axis of the moving frame. We get the coordinates $x$, $y$ of the point $P$ in the fixed system by
$$
x=p+u\, \cos \varphi -v\, \sin \varphi\,, \qquad y=q+u\, \sin \varphi +v\, \cos \varphi.
$$
Here $p,q,\varphi$ are functions of a quantity $t$ which determine the motion. (For example, $t$ can denote the time, or any other metric parameter.) Assume that $\varphi (t)$ is not zero on an interval of $t$. Then it can be inverted, and $p,q$ can also be considered as a function of $\varphi$. (This assumption says that our motion cannot contain translations in that domain. We call such a motion \emph{non-translative planar motion}.) The derivative of the coordinate functions with respect to $\varphi$ gives the coordinates of the velocity vector of the point $P$.
It is more convenient to use vector equality, and hence we introduce some further notion. Let
$$
\R(\varphi)=\left(\begin{array}{cc}
                \cos \varphi & -\sin \varphi \\
                \sin \varphi & \cos \varphi
                \end{array}
            \right)
$$
denote the rotation about the origin with signed angle $\varphi$. Then the first equation array has the form
$$
{\bf x}={\bf p}+\R(\varphi){\bf u}
$$
If $\Q=\R(\pi/2)$ denotes the rotation with $\pi/2$, we have the following rules:
$$
\Q^2=-\E, \quad \Q^3=\Q^{-1}=\overline{\Q}=-\Q, \quad \Q^4=\E,
$$
where $\E$ is the unit matrix. We denote by dot the \emph{derivative with respect to} $\varphi$, which means in this section the Euclidean arc-length parameter. It is clear that
$$
\dot{\R}=\Q\R, \quad \dot{(\R^{-1})}=-\Q\R.
$$
For every value of $\varphi$ there is precisely one point ${\bf u}_0$ of the moving plane for which the velocity vector vanishes. This is
$$
{\bf u}_0=\Q\R^{-1}\dot{\bf p}.
$$
This point ${\bf u}_0$ of the moving plane is a so-called \emph{instantaneous center ({\rm or} instantaneous pole)} of the motion, and the set of these points is the \emph{moving polode}, or curve $\gamma'$ of instantaneous poles, of the moving plane. The points of the moving polode can also be obtained in the frame as rest. These points ${\bf x}_0$ are of the form
$$
{\bf x}_0={\bf p}+\R{\bf u}_0={\bf p}+\Q\dot{{\bf p}},
$$
and form the so-called  \emph{fixed polode}, or curve $\gamma $ of instantaneous centers, in the fixed plane. We examine the motion with respect to the point ${\bf x}_0$. If ${\bf x}$ is arbitrary, then ${\bf x}-{\bf x}_0=\R{\bf u}-{\bf Q}\dot{\bf p}$, and using the equality  $\dot{\bf x}=\dot{\bf p}+\Q\R{\bf u}$, we have $\Q \dot{\bf x}=\Q \dot{\bf p}+\Q\R{\bf u}$. Since ${\bf x}-{\bf x}_0=\R {\bf u}-\Q\dot{\bf p}$, we get that
$$
\dot{\bf x}=\Q({\bf x}-{\bf x}_0).
$$
Hence the velocity vector of the motion at the point ${\bf x}$ is orthogonal to the position vector from ${\bf x}_0$ to ${\bf x}$. This implies that the moving system in the given moment is a rotation about the center ${\bf x}_0$. Observe that the velocity vectors of the two polodes at their common point agree; in fact,
$$
\dot{{\bf u}}_0=\dot{\Q\R^{-1}\dot{\bf p}}=\R^{-1}\dot{\bf p}+\Q\R^{-1}\ddot{\bf p}=\dot{\bf x}_0.
$$
Hence the arc-length elements of the two curves agree and we get that in every moment the two curves are touching, and their arc-lengths calculated from a point $\varphi_0$ to the point $\varphi $ have the same value. Hence the moving polode $\gamma'$ \emph{rolls without slipping} (or without friction) on the fixed polode $\gamma$, and this is the only rolling process which corresponds to the given motion of the planes. Hence we got the fact \emph{that every non-translatory planar motion of a rigid mechanical system in the plane can be considered as the rolling process of a curve rigidly connected with the system on a fixed curve in the plane}.
This motivates the so-called main theorem of planar Kinematics, namely: At every moment, any constrained non-translatory planar motion can be approximated (up to the first derivative) by an \emph{instantaneous rotation}. The
center of this rotation is called the \emph{instantaneous pole}. Thus, for each
position of the moving plane, we generally have exactly one point with velocity zero (as a
result of that, the instantaneous pole is also called velocity center).

This theorem leads to an interesting class of curves in the Euclidean plane.
Given a curve $\gamma'$ attached to a plane $\Sigma'$ which is moving so that the curve rolls, without slipping, along a given curve $\gamma $ attached to a fixed plane $\Sigma$ occupying the same space, then a point $P$ attached to $\Sigma'$ describes a curve in $\Sigma $ called a \emph{roulette}.

Based on this rolling process we can rewrite the definition of the motion of rigid systems. Observe that every planar motion implies the motion of all points of the moving plane with respect to the fixed one. These orbits are said to be roulettes. Thus, for the studied motion we consider two curves, also called \emph{polodes}, and a suitable rolling process to determine the motion of a singular point. For this purpose a method is needed to determine the fixed position of the point $P$ with respect to the moving polode. A usual method is to give a line through the point $P$ which intersects the moving polode in the point $Q$ and fixes the distance of $P$ and $Q$ and the angle of the line $PQ$ with the tangent line $t_Q$ of the moving polode at $Q$. Hence the choice of $Q$ on the moving polode is arbitrary. Fix $Q={\bf w}(0)$ and $P={\bf x}(0)$.  The points of the roulette ${\bf w}(s)$ of $Q$ can be obtained by the composition of the following transformations: translate the point $\gamma'(s)$ into the origin, rotate the image of the point of $\gamma(0)$ about the origin by the angle $\varphi(s)=\left( \dot{\gamma}(s),\dot{\gamma'}(s)\right)\angle $, and translate the obtained point by $\gamma(s)$. Hence the equality of the roulette of $Q$ in the fixed system is
$$
{\bf w}(s)=\R(\varphi(s))(-\gamma'(s))+\gamma(s)=\gamma(s)-\R(\varphi(s))(\gamma'(s)).
$$
Since the roulette ${\bf x}(s)$ of the point $P$ can be described by the formula ${\bf x}(s)={\bf w}(s)+\R(\varphi(s)){\bf p}$, we get
$$
{\bf x}(s)=\gamma(s)+\R(\varphi(s))\left({\bf p}-\gamma'(s)\right).
$$
This means that if we have two touching arcs $\gamma (s)$ and $\gamma'(s)$  of a plane $\Sigma $, and we associate to the second arc a moving plane $\Sigma '$ in which its position is fixed, then the rolling process of $\gamma'(s)$ on $\gamma(s)$ (locally) uniquely determines an orbit of every point of $\Sigma'$. In the Euclidean plane equation the above equation shows that in every moment with respect to varying $p$ we have an isometry. Hence the rolling process of the arcs determines a rigid motion of the plane $\Sigma'$. This representation is locally unique, since a rigid motion uniquely determines its polodes. Hence we have if $\gamma,\gamma':[0,\beta]\rightarrow \mathbb{R}^2$ are two simple Jordan arcs with common touching point $\gamma(0)=\gamma'(0)$ such that $s$ is the arc-length parameter of both of them (considered from the points $\gamma(0), \gamma'(0)$ to the point $\gamma(s)$, $\gamma'(s)$, respectively), then for every $s\in [0,\beta]$ we have an isometry $\Phi_s$ sending the original position vector ${\bf p}$ into the instantaneously position $\Phi_s(p)$. If $\gamma$ and $\gamma'$  have, for all $s\in [0,\beta]$, unique tangents at their points $\gamma(s)$ and $\gamma'(s)$, respectively, then, for all $s\in[0,\beta]$, $\Phi_s$ is uniquely determined and can be described by the vector equation
$$
\Phi_s({\bf p})=\gamma(s)+\R(\left( \dot{\gamma}(s),\dot{\gamma'}(s)\right)\angle)\left({\bf p}-\gamma'(s)\right)\,.
$$
Here $\dot{\gamma}(s)$ and $\dot{\gamma'}(s)$ denote the unit tangent vectors at the point $\gamma(s)$ and $\gamma'(s)$, respectively, and $\R(\theta)$ is the rotation with the angle $\theta$. For fixed ${\bf p}$, the graph of the function $\Phi_{(\cdot)}({\bf p}):[0,\beta]\rightarrow \Sigma $ is said to be the \emph{roulette} of the point $P={\bf p}\in \Sigma$ for the rigid motion given by the system of isometries $\{\Phi_s:\, s\in[0,\beta]\}$.

The most important results in this theory are the so-called Euler-Savary equations, which compare the curvatures of the moving polode, fixed polode and the corresponding roulettes. We use the following forms of them (compare with Fig. \ref{EulerSavaryfig}):
\begin{figure}
 \centering
   \includegraphics[]{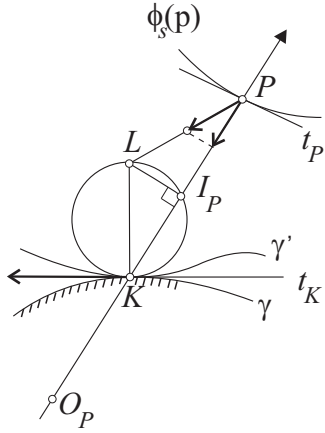}
 \caption{Roulette}
\label{EulerSavaryfig}
\end{figure}
The \emph{first Euler-Savary equation} is
$$
 \frac{1}{r'}-\frac{1}{r^\star}=\frac{1}{\alpha \sin \nu},
$$
where $r',r^\star$ are the curvature radius of the fixed polode at the instantaneous center and the curvature radius of the roulette at the examined point, respectively, $\alpha$ is the length of the velocity vector and $\nu$ is the angular velocity of the motion.
The \emph{second Euler-Savary equation} says that
$$
\frac{1}{r'}-\frac{1}{r''}=\frac{1}{\alpha},
$$
where $r''$ is the curvature radius of the moving polode at the instantaneous center. From these equations we can get a common form of the two equations which is
$$
\frac{1}{r'}-\frac{1}{r''}=\left(\frac{1}{r'}-\frac{1}{r^\star}\right)\sin \nu.
$$

\subsection{Roulettes in spherical geometry}

Spherical geometry of dimension $2$ is the geometry of a $2$-sphere embedded into a $3$-dimensional Euclidean space. The motions of the sphere can be get from the three-dimensional special orthogonal linear group $\mathrm{SO}(3)$, containing those orthogonal transformations which hold the orientation. Algebraically, this subgroup contains those linear transformations of the $3$-dimensional Euclidean space which determinant has value $1$. Euler observed that these transformations always has an eigenvector with eigenvalue $1$ showing that a motion of the sphere is a rotation about a line which is called by the rotational axis of the transformation. Hence the rigid motions of the sphere not only the sequence of certain instantaneous rotations but the sequence of proper rotations of the sphere. In fact, in Euclidean geometry we should combine the rotation of the coordinate system with the translation of its origin, contrary on the sphere the translation part is also a rotation (at another point as the fixed origin) hence the needed composition gives the product of three rotations which itself also a rotation (at a new origin). To characterize the motion we have to give in a moment a rotation, hence we have a function of the time, mapping the domain time-interval into $\mathrm{SO}(3)$, continuously.

There are several papers on spherical roulettes because spherical kinematics is an intermediate step between planar and spatial kinematics. The earliest one is found by me, the work written by Garnier in 1956 \cite{garnier}. Considering one and two parameters spherical motions in
Euclidean space, Muller gave the relations for absolute, sliding, relative velocities and pole curves of these motions. In addition
to that he expressed the corresponding Euler–Savary formula related to the trajectory curves of these $1$-parameter spherical motions \cite{muller}. I would like to mention here also the paper of Chiang \cite{chiang} which contains the results of the most important theorems using spherical kinematics, especially the spherical form of Euler-Savary equations. For the verification of the results Chiang proposed the books \cite{dittrich} and \cite{chiang1}. I would like to mention two nice recent papers. The first on by M.A. Gungor, S. Ersoy, M. Tosun \cite{gungor} containing more information on the history, too. The second one written by Brunnthaler, Husty and Schr\"ocker on four-bar mechanism (\cite{brunnthaler-husty-schroecker}).

\begin{figure}[ht]
 \centering
   \includegraphics[scale=0.5]{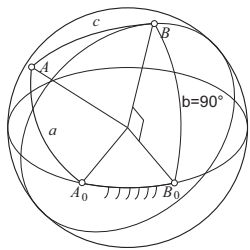}
 \caption{Spherical rotary slider-crank}
\label{chiangfig1}
\end{figure}

In the spherical plane the concept of roulettes also strongly connected to the bar-joint kinematics, in Chiang's paper we can find the complete lists of motions induced by bar-joint frameworks. As Chiang said:
\emph{Spherical four-bar mechanisms are similar to planar four-bar mechanisms. There are also spherical crank-rockers, drag-links and double-rockers, and even also spherical 'slider-cranks'.
However, as there is no translation motion in spherical kinematics, because all spherical motions are rotations, there exists therefore no real spherical 'slider'. What we call a 'slider' here exists simply because a joint on the slider is moving along a great circle arc. In fact this joint is $90^\circ$ apart from its axis of rotation. It is interesting to note the spherical rotary slider-crank as that shown
in Fig.\ref{chiangfig1}, and there is no planar counterpart.}

\begin{figure}[ht]
 \centering
   \includegraphics[scale=0.5]{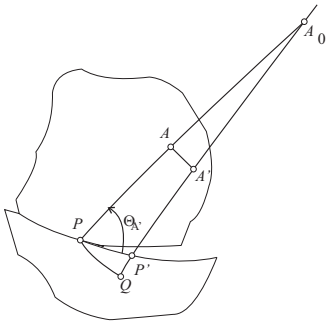}
 \caption{Spherical Euler-Savary equation}
\label{chiangfig2}
\end{figure}

Using the notation of Chiang we can say the following. Assume that the polodes contact at $P$. Denote by $A$ the point of the roulette, by $\Theta_{A}$ the spherical argument of the spherical polar coordinates of $A$, with the pole-tangent great circle as the spherical 'polar line' and by $\Theta=-u/\omega$ the quantity corresponding to the diameter of the inflection circle in planar kinematics. (Note that $u$ is the pole changing velocity, and $\omega$ is the angular velocity (magnitude) of the moving body, respectively.) Finally, if $A_0$ means the center of curvature of the roulette then we have the following (geometric form) of the first Euler-Savary equation (see in Fig \ref{chiangfig2}):
\begin{theorem}[\cite{chiang}]
By the notation introduced above we get the following connection among the spherical lengths in question:
$$
\frac{1}{\tan PA}-\frac{1}{\tan PA_0}=\frac{1}{\Theta\sin\Theta_A}
$$
\end{theorem}

\subsection{Roulettes in the hyperbolic plane}

It is very surprising but I can not find any paper on the internet on hyperbolic roulettes\footnote{This does not means that there is no results on roulettes in hyperbolic geometry. Thank you for Hans-Peter Schr\"ocker who recommended some works can be cited in this section.}. It seems to be that the spherical approach can works also in the Lorentzian space (using it for its imaginary unit sphere). This analytically can solves the case of the hyperbolic plane. Since roulettes can be defined by non-translatory motions, the basic problem of the general transition to hyperbolic geometry (that is we have two types of isometries can work as a translation) will not be occur. We can predict also the form of the Euler-Savary equation knowing the spherical one of the preceding paragraph.  Using the connections between the metric of the two planes mentioned in the introduction and the notation of  Chiang's corresponding to the motion we think that the required equation is:
$$
\frac{1}{\tanh PA'}-\frac{1}{\tanh PA'_0}=\frac{1}{\Theta\sin \Theta_A'}.
$$
(So the polodes contact at $P$, $A'$ is the point of the roulette $\Theta_A' $ is the hyperbolic argument of the hyperbolic polar coordinates of $A'$, with the tangent line and $\Theta=-u/\omega$ is a quantity corresponding to the diameter of the inflection circle in planar kinematics, $u$ is the pole changing velocity, and $\omega$ is the angular velocity (magnitude) of the moving body.

For the interested reader we propose the works \cite{garnier}, \cite{frank}, \cite{tolke 1}, \cite{tolke 2}, \cite{tolke 3} on hyperbolic kinematics.

\subsection{Roulettes in the Lorentzian-plane}

Ergin \cite{ergin} considering the Lorentzian plane instead of the Euclidean plane, and introduced the one-parameter planar motion in the Lorentzian
plane and also gave the relations between both the velocities and accelerations. Y\"uce and Kuruoglu in \cite{yuce} using hyperbolic numbers reproduce the results of Ergin and in analogy with complex motions as given by M\"uller \cite{muller}, defined one parameter motions in the Lorentzian plane. They calling is hyperbolic plane is not good, because of hyperbolic numbers can be identified with the Lorentzian plane and not with the hyperbolic plane. The relations between absolute, relative, sliding velocities (and accelerations) and pole curves was discussed, too. In the Lorentzian plane Euler-Savary formula is given in references, \cite{ergut} and \cite{ikawa}. We now refer the Theorem 4.1 from Ikawa's paper. By $c(t) = (x(t), y(t))$, Ikawa denoted the orthogonal coordinate representation of a curve $c(t)$. The vector field
$\frac{\mathrm{d}c(t)}{\mathrm{d}t}$ is called the tangent vector field of the curve $c(t)$. If the tangent vector field $X$ of $c(t)$ is a spacelike, timelike, or null, then the curve $c(t)$ is called spacelike, timelike, or null, respectively. Ikawa did not consider null curves and proved the following Euler-Savary types theorem:
\begin{theorem}[\cite{ikawa}]
On the Lorentzian plane $L^2$, suppose that a curve $c_R$ rolls without slipping along a curve $c_B$. Let $c_L$ be a locus of a point $P$ that is relative to $c_R$. Let $Q$ be a point on $c_L$ and $R$ a point of contact of $c_B$ and $c_R$ corresponds to $Q$ relative to
the rolling relation. By $(r;\varphi)$, we denote a polar coordinate of $Q$ with respect to the origin $R$ and the base line ${c'_B}|_{R}$. Then curvatures $k_B$; $k_R$ and $k_L$ of $c_B$; $c_R$ and $c_L$, respectively, satisﬁes
$$
rk_L = \pm 1-\frac{\cosh \varphi}{r|k_B-k_R|} \mbox{ when $c_L$ is space-like;}
$$
$$
rk_L = \pm 1+\frac{\cosh \varphi}{r|k_B-k_R|} \mbox{ when $c_L$ is time-like.}
$$
\end{theorem}

\subsection{Roulettes in Minkowski normed plane}

The Euler-Savary equation clarifies the relation between the curvatures of the fixed and the moving polodes (and the respective roulettes) of a rolling process without friction, determining a planar motion of a rigid body (and vice versa). This connection follows from Euler's theorem of classical mechanics which states that every planar motion can be considered as the sequence of instantaneous rotations, whose centers give the fixed polode. Thus, the Euler-Savary equation has been investigated also in physics, in various contexts (Lorentz spatial motions, elliptical harmonic motions, homothetic motions in the complex plane, etc.). However, one reason that there are no deeper and wider investigations on Euler-Savary applications in physics might be the fact that there are no rotations in the classical meaning. From a certain point of view, in the paper \cite{balestro-gha-martini} filled this gap by proposing a concept of rotation which mathematically describes a planar motion with respect to normed planes. Thinking about the concept of motion in a wider context, it can be allowed that the duration of the motion of a body a little bit changes its shape. This "changing of shape" is determined by the concrete moving (and thus by the participating polodes and roulettes). In order to define a concept of rotation for a Minkowski plane, the authors started with extending the definition of Brass by considering Borel measures in a larger class of curves, not only in the unit circle, and derived  angle measures for normed planes from it. (The recent survey \cite{balestro-gha-martini-teixeira} contains a lot of information on angle measures in Minkowski plane.)

Let $\gamma \subseteq X$ be a closed Jordan curve which is starlike with respect to a point $p$ of the interior of the region bounded by $\gamma$. An \emph{angle measure with respect to such a Jordan curve} is a (normalized) Borel measure $\mu_{\gamma}$ on $\gamma$ for which the following properties hold:

\noindent\textbf{(a)} $\mu_{\gamma}(\gamma) = 2\pi$;

\noindent\textbf{(b)} for any $q \in \gamma$ we have $\mu_{\gamma}(\{q\}) = 0$; and

\noindent\textbf{(c)} any non-degenerate arc of $\gamma$ has positive measure.

\begin{figure}[ht]
\centering
\includegraphics[height=3cm]{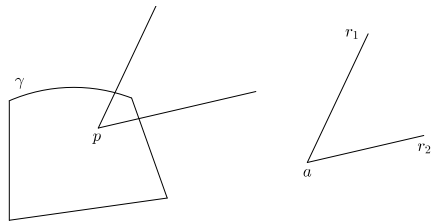}
\caption{The generalized angle measure given by $\mu_{\gamma}$ and $p$}
\label{figgenang}
\end{figure}

An angle measure defined in this way provides a translation invariant measure of \emph{angles} in the plane, which can be defined as the convex hulls of two rays with the same starting point. Given an angle $(r_1,r_2)\angle$ with apex $a$, it can be defined its \textit{generalized angle measure} $\mu_{\gamma,p}(r_1,r_2)$ to be the measure $\mu_{\gamma}$ of the arc determined on $\gamma$ by the image of $(r_1,r_2)\angle$ via the translation $x \mapsto x - a + p$. Figure \ref{figgenang} illustrates this concept.

Using this notion of generalized angle measure we define now the generalized rotations in Minkowski planes.

\begin{defi}\label{generalrotation} Let $(X,||\cdot||)$ be a Minkowski plane and let $\gamma$ be a closed Jordan curve which is starlike with respect to a point $p$ of the interior of the region bounded by $\gamma$. Let $\mu_{\gamma,p}$ be a generalized angle measure as in the previous definition. A \textit{general rotation} (with respect to $\mu_{\gamma,p}$) is a transform $\mathrm{rot}_{\mu_{\gamma,p}}:X \rightarrow X$ for which the following three properties hold:

\normalfont
\noindent\textbf{(a)} The transform $\mathrm{rot}_{\mu_{\gamma,p}}$ leaves invariant the pencil $\mathcal{R}(p)$ of rays with origin in $p$. In other words, if $r \subseteq X$ is a ray with origin $p$, then $\mathrm{rot}_{\mu_{\gamma,p}}(r)$ is also a ray with origin $p$.

\noindent\textbf{(b)} For each $\alpha > 0$, $\mathrm{rot}_{\mu_{\gamma,p}}$ leaves invariant the homothetic curve $\gamma_{\alpha,p} := p + \alpha(\gamma - p)$, i.e., for such a curve we have $\mathrm{rot}_{\mu_{\gamma,p}}\left(\gamma_{\alpha,p}\right) \subseteq \gamma_{\alpha,p}$.

\noindent\textbf{(c)} The function $r \in \mathcal{R}(p) \mapsto \mu_{\gamma,p}\left(\mathrm{rot}_{\mu_{\gamma,p}}(r),r\right)$ is constant. Intuitively, $\mathrm{rot}_{\mu_{\gamma,p}}$ ``rotates every ray of $\mathcal{R}(p)$ by a same angle".
\end{defi}

Notice that a general rotation can be considered as acting in the space of directions of $X$. Indeed, the set $\mathcal{R}(p)$ can be seen as this space. For a class $\mathcal{R}(\gamma,\mu,p)$ we have the following properties:
\begin{itemize}
\item Regarding composition, $\mathcal{R}(\gamma,\mu,p)$ is an abelian group. More precisely, we have $\mathrm{rot}_{\theta_1}\circ\mathrm{rot}_{\theta_2} = \mathrm{rot}_{\theta_1\oplus \theta_2}$, where $\oplus$ is the sum modulo $2\pi$.

\item For any $q\in\gamma$, the application $l\mapsto\mathrm{rot}_{\theta}(q)$ is a bijection from $[0,2\pi)$ to $\gamma$.
\end{itemize}

We highlight an interesting fact: The standard Euclidean rotation group can be obtained in any Minkowski plane. We just have to consider the group $\mathcal{R}(\gamma,\mu,o)$ where $\gamma$ is the \emph{L\"{o}wner ellipse}, which is defined as the ellipse of maximal volume contained in $B$, and $\mu$ is the measure given by twice the area of its sectors.

In \cite{balestro-gha-martini} there are two examples of general rotations in the Euclidean plane. The first one relies on an area-based measure for an ellipse, which is clearly well defined. In the second used the arc-length measure referring to a nephroid. Here we quote only the first one.

\begin{figure}[ht]
\centering
\includegraphics[]{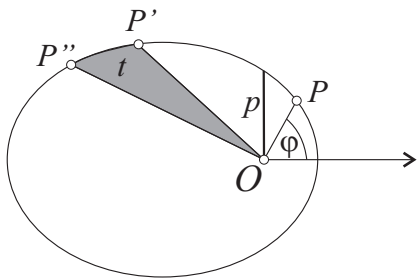}
\caption{Area-based rotation and the Kepler's model}
\label{heliofig}
\end{figure}

Consider the Euclidean plane and the system of ellipses with common focus at the origin $O$ and with major axis on the $x$-axis of the coordinate system, such that the positive half-line of $x$ contains the closest point of the ellipse (see Fig. \ref{heliofig}). In that polar coordinate system (which is called the heliocentric coordinate system for the ellipse), for which the ray $\varphi=0$ is the positive half axis $x$, we can write the radial function $r(\varphi)$ of the ellipse $G$ by the formula
$$
r(\varphi)=\frac{p}{1+\varepsilon\, \cos \varphi},
$$
where $p$ is the semi-latus rectum of the ellipse and $\varepsilon$ is the eccentricity of it, respectively. Let $\mu((\varphi',\varphi'')\angle)$ be the area of the sector enclosed by $\varphi'$, $\varphi''$, and $G$ be the arc between these lines. Hence
$$
\mu((\varphi',\varphi'')\angle)=\frac{1}{2}\int\limits_{\varphi'}^{\varphi''}\left(\frac{p}{1+\varepsilon\, \cos \varphi}\right)^2\mathrm{d}\varphi.
$$

With respect to $\mu$ and $G$ from above, for every real number $0\leq t\leq 2\pi$ there is a generalized rotation of the Euclidean plane about $O$ with this angle $t$. By Kepler's second law about planetary motions, the angle $t$ of a generalized rotation is proportional to the time of the motion of the planet. Hence the generalized rotation with angle $t$ maps the current position $P'$ of the planet to that point $P''$ of the orbit where the planet arrives after time $t$.

The principle of measuring the angle proportional to the area of the sector intersected by the angle domain from the basic disk $\left(G\cup \mathrm{int} G\right)$ works in all Minkowski planes and for all basic curves $G$. Note that in the Euclidean plane with the unit circle as basic curve, this choice of $\mu$ gives the usual angle measure, and that we get the usual rotations as generalized rotations by choosing $P$ to be the origin $O$. An advantage of this choice is affine invariance, but there is also a big disadvantage. Namely, the length of the arc $G$ containing the domain of the angle cannot be calculated easily from this angle measure. (As a known example, we note that the calculation of the arc-length of an ellipse leads to a complete elliptic integral of second kind, which has no closed-form solution in terms of elementary functions.)

If we consider two curves $\gamma$ and $\gamma'$, then we have to use a  suitable lower subscript for the curvature function. We also have the concept of \emph{curvature radius} $r_{\gamma}$ which is, as well-known, the reciprocal value of the curvature at the given point $K=\gamma(s)$. With these notions we are able to formulate

If the unit circle of the Minkowski plane is two times continuously differentiable, then the following equality holds:
\begin{equation}\label{equ:Euler-Savary}
\chi_{\gamma}-\chi_{\gamma'}=\frac{1}{r_{\gamma}}-\frac{1}{r_{\gamma'}}=\frac{\sigma(T_K)}{\sigma^2(t_K)}\frac{1}{\alpha_K}\,.
\end{equation}
Here $r_{\gamma}$ is the curvature radius of the fixed polode at its point $K=\gamma_s$, $r_{\gamma'}$ is the curvature radius of the moving polode at its point $K=\gamma'_s$, and $\alpha_K$ is the length of the common velocity vector of the fixed and moving polodes at the moment $s$ and at the instantaneous pole $K=\gamma(s)=\gamma'(s)$.

Deeper investigation leads to the following geometric form of the first Euler-Savary theorem (see Fig \ref{EulerSavaryfig}).

The instantaneous center $K$ and the curvature center $O_P$ of the roulette at its point $P\ne K$ satisfy the equality
$$
\|\overrightarrow{O_PP}\|=\frac{\|\overrightarrow{KP}\|^2}{\|\overrightarrow{I_PP}\|},
$$
where the second intersection point of the path normal line at $P$ with the inflection curve is the point $I_P$.

This yields the \emph{combined formula of the two Euler-Savary equations}, namely
\begin{theorem}{\cite{balestro-gha-martini}}
With the above notation we have
$$
\left(\frac{1}{KP}-\frac{1}{KO_P}\right)\mathrm{sm}(g(K,P),t_K)\frac{\sigma^2(g(K,P))}{\sigma^2(t_K)\sigma(g(K,L))}= \dot{\varphi}(0)\left(\chi_{\gamma}-\chi_{\gamma'}\right)=\frac{\dot{\varphi}(0)}{\sigma^2(t_K)}\,\frac{1}{\alpha_K},
$$
where we assume that $\sigma(T_K)=\mathrm{area}B=1$ and $\mathrm{sm}$ means the sine function of Busemann.
\end{theorem}

\end{document}